\documentclass[11pt]{article}
\usepackage{amsmath}
\usepackage{amssymb}
\usepackage{amsfonts}

\usepackage{epsf}

\usepackage[OT2,OT1]{fontenc}
\newcommand\cyr{%
\renewcommand\rmdefault{wncyr}%
\renewcommand\sfdefault{wncyss}%
\renewcommand\encodingdefault{OT2}%
\normalfont
\selectfont}
\DeclareTextFontCommand{\textcyr}{\cyr}

\begin{document}

\def\f#1#2{\mbox{${\textstyle\frac{#1}{#2}}$}}

\title{Zeroes of the Swallowtail Integral}
\date{}
\author{David Kaminski}
\maketitle

\abstract{The swallowtail integral 
\begin{displaymath}
S(x,y,z) = \int_{-\infty}^{\infty} \exp[i(u^5 + xu^3 + yu^2 + zu)] \, du
\end{displaymath}
is one of the so-called canonical diffraction integrals used in optics, and
plays a role in the uniform asymptotics of integrals exhibiting a confluence
of up to four saddle points. In a 1984 paper by Connor, Curtis and Farrelly, the authors make a number of remarkable observations regarding the zeroes of
$S(x,y,z)$, including that its zeroes occur on lines in $xyz$-space, and
that the zeroes of $S(0,y,z)$ lie along the line $y = 0$. These assertions
are based on numerical evidence and the asymptotics of $S(0,0,z)$. We examine these assertions more completely and provide additional
detail on the structure of the zeroes of $S(x,y,z)$.}

\section{Introduction}
The swallowtail integral
\begin{equation}\label{sw}
S(x,y,z) = \int_{-\infty}^\infty \exp[i(u^5 + xu^3 + yu^2 + zu)]\, du
\end{equation}
is one of the suite of integrals used in asymptotics for the construction of uniform asymptotic expansions of integrals in which four saddle points coalesce \cite[Ch.~VII, \S 6]{W}, and has a home in optics where the integral is used to describe diffraction phenomena \cite{BK}. These types of integrals were the object of considerable study in the 1980s and 1990s, but still surface in the literature on occasion in more recent years (see, for example \cite{CH} and \cite{Nye}), and have even found a home in the {\em NIST Handbook of Mathematical Functions} \cite[Ch.~36]{NIST}.

When working with $S(x,y,z)$ it is often less notationally cumbersome to use a slightly rescaled version,
\begin{equation}\label{Q}
Q(x,y,z) = \int_{-\infty}^\infty \exp[i(\f{1}{5}t^5 + \f{1}{3}xt^3 + \f{1}{2}yt^2 + zt)]\, dt.
\end{equation}
$S$ and $Q$ are easily related by
\begin{equation}\label{SQ}
   S(x,y,z) = \frac{1}{5^{1/5}}Q(\frac{3x}{5^{3/5}}, \frac{2y}{5^{2/5}},
   \frac{z}{5^{1/5}})
\end{equation}
and
\begin{equation}\label{QS}
   Q(x,y,z) = 5^{1/5}S(\frac{5^{3/5}x}{3}, \frac{5^{2/5}y}{2}, 5^{1/5}z).
\end{equation}
Both $S$ and $Q$ enjoy the symmetry 
\begin{equation}\label{symmetry}
   S(x,-y,z) = \overline{S(x,y,z)}, \qquad Q(x,-y,z) = \overline{Q(x,y,z)},
\end{equation}
where the overline represents complex conjugation,
so all our work involving $S$ and $Q$ may be restricted to $y \geq 0$ without loss of generality.
Nye \cite{Nye} uses ${\displaystyle \frac{1}{\sqrt{2\pi}}Q(x,y,z)}$ as his swallowtail integral.

We note that through an application of Jordan's lemma, the integration contour can be deformed into one that begins at ${\infty}e^{9{\pi}i/10}$ and ends at ${\infty}e^{{\pi}i/10}$, so that the integrand undergoes exponential decay as $|t| \rightarrow \infty$ along the integration contour.

Of interest to us in the present work is the distribution of zeroes of $Q$ for large values of the parameters $x$, $y$ and $z$, which we take to be real. In particular, our attention will initially focus on the remarkable observations made in \cite{CCF} that the zeroes of $S$ (and therefore of $Q$) lie along lines in $xyz$-space, and that the zeroes of $S(0,y,z)$ are to be found along the $z$-axis in the $yz$-plane. The setting in \cite{CCF} is one of numerical computation, and the authors appear to be making the claim based on computated values of $S$ they made, and draw attention to the connexion of the asymptotics of $S(0,0,z)$ and location of the computed zeroes.

Of related interest is a set of observations made in \cite{Nye} where reference is made to saddle points of $Q(x,0,z)$. Some additional commentary on Nye's work is provided at the close of this paper.

\section{$Q(0,0,z)$}
With 
\begin{equation}\label{fpm}
   f_{\pm}(t) = \f{1}{5}t^5 \pm t
\end{equation} 
we see that a change of variable $t \mapsto |z|^{1/4}t$ permits us to 
write
\begin{equation}\label{Spm}
   Q(0,0,\pm|z|) = |z|^{1/4}\int_{\infty e^{9\pi i/10}}^{\infty e^{\pi i/10}}
   e^{i |z|^{5/4}f_{\pm}(t)} dt.
\end{equation}
It follows that for the case of $z > 0$, the saddle points for $Q(0,0,z)$
are $t = \pm e^{\pm\pi i/4}$, or
\begin{equation}\label{tkp}
   t_k = i^i e^{\pi i/4}, \qquad k = 0, 1, 2, 3.
\end{equation}
For $z < 0$, the corresponding saddles for $Q(0,0,-|z|)$ are $t = \pm 1, \pm i$, or
\begin{equation}\label{tkm}
   t_k = i^k, \qquad k = 0, 1, 2, 3.
\end{equation}
We observe that
\begin{displaymath}
   f_{+}(t_k) = \f{4}{5}i^k e^{\pi i/4} \quad
   \textrm{and} \quad
   f_{-}(t_k) = -\f{4}{5}i^k
\end{displaymath}
where the selection of $t_k$ is made in (\ref{tkp}) or (\ref{tkm}) according
to the sign of $z$.

The steepest descent paths through these saddles in each case ($z > 0$ or
$z < 0$) are among the steepest curves through the saddles. For $z > 0$,
\begin{eqnarray*}
   i\bigl(f_{+}(t) - f_{+}(t_k)\bigr) & = &
   2i^{3k+1}e^{3\pi i/4}(t - t_k)^2 - 2(-1)^k(t-t_k)^3 \\
   & & \quad + i^{k+1}e^{\pi i/4}(t-t_k)^4 + \f{1}{5}i(t-t_k)^5,
\end{eqnarray*}
and for $z < 0$,
\begin{eqnarray*}
   i\bigl(f_{-}(t) - f_{-}(t_k)\bigr) & = &
   2i^{3k+1}(t-t_k)^2 + 2i^{2k+1}(t-t_k)^3 \\
   & & \quad + i^{k+1}(t-t_k)^4+\f{1}{5}i(t-t_k)^5,
\end{eqnarray*}
where, again, the selection of saddles $t_k$ is made from (\ref{tkp}) or (\ref{tkm}) according as $z > 0$ or $z < 0$. The steepest curves for $z > 0$ passing through $t_0$ and $t_1$ are depicted in Fig.~1 (the origin in each illustration in the figure corresponds the saddle point, so the left illustration in Fig.~1 depicts $t_1$ at the origin, and the right illustration in Fig.~1 has $t_0$ located at the origin). 

For $z > 0$, the steepest descent curve through $t_1$ is that curve $\Gamma_1^{+}$  beginning at $\infty e^{9\pi i/10}$ which passes through $t_1$ and then rises to $\infty i$; the steepest descent curve through $t_0$ is that curve $\Gamma_0^{+}$ beginning at $\infty i$ which passes through $t_0$ and then ends at $\infty e^{\pi i/10}$. Thus, for 
$S(0,0,z)$ with $z > 0$, the integration contour can be deformed into the sum of these two steepest descent curves, $\Gamma_1^{+} + \Gamma_0^{+}$, yielding two relevant saddle points for the asymptotics of $S$.

\begin{figure}[h]%
\begin{center}
\leavevmode
\epsfxsize=2in\epsfbox{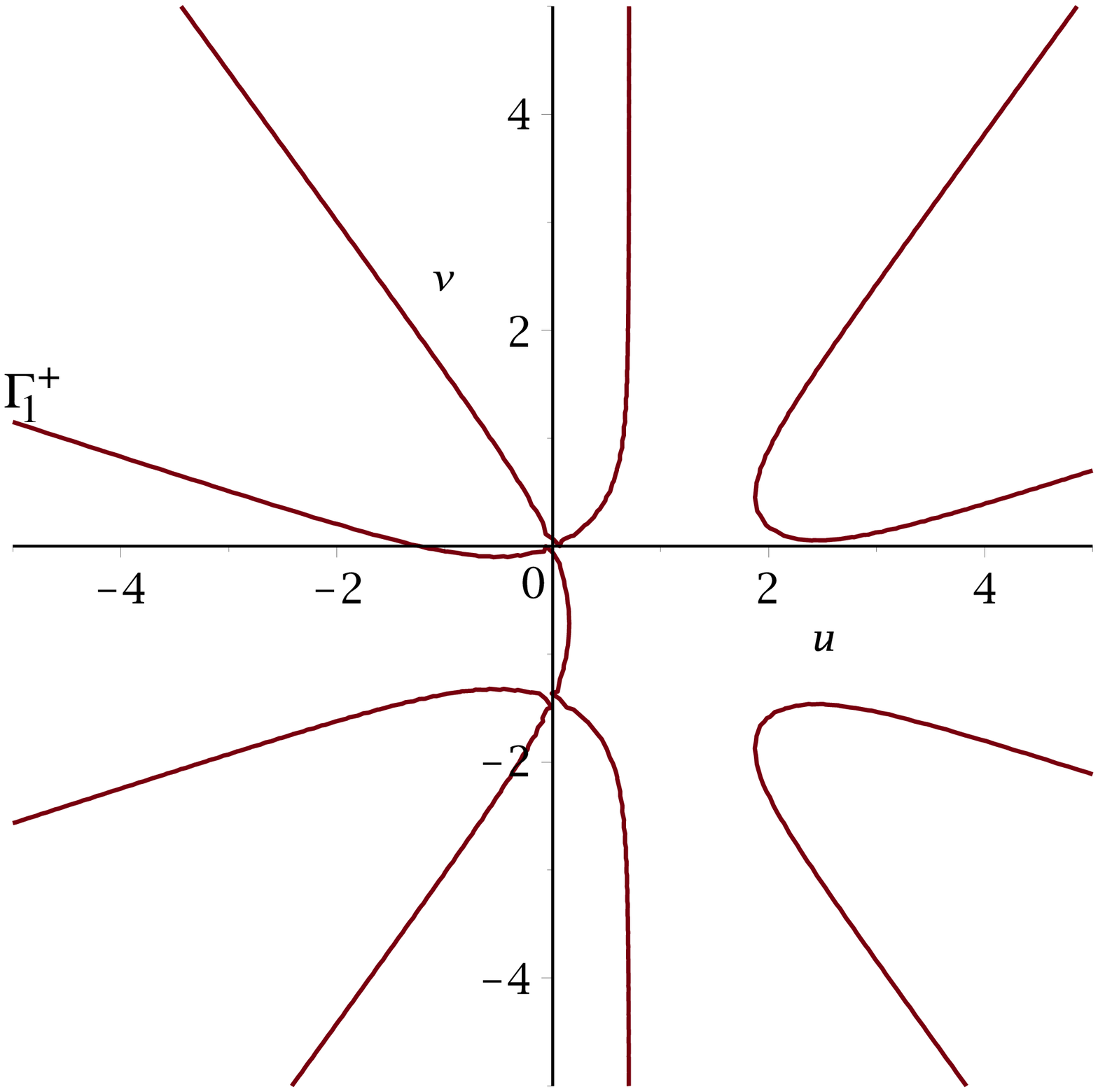}\quad\epsfxsize=2in\epsfbox{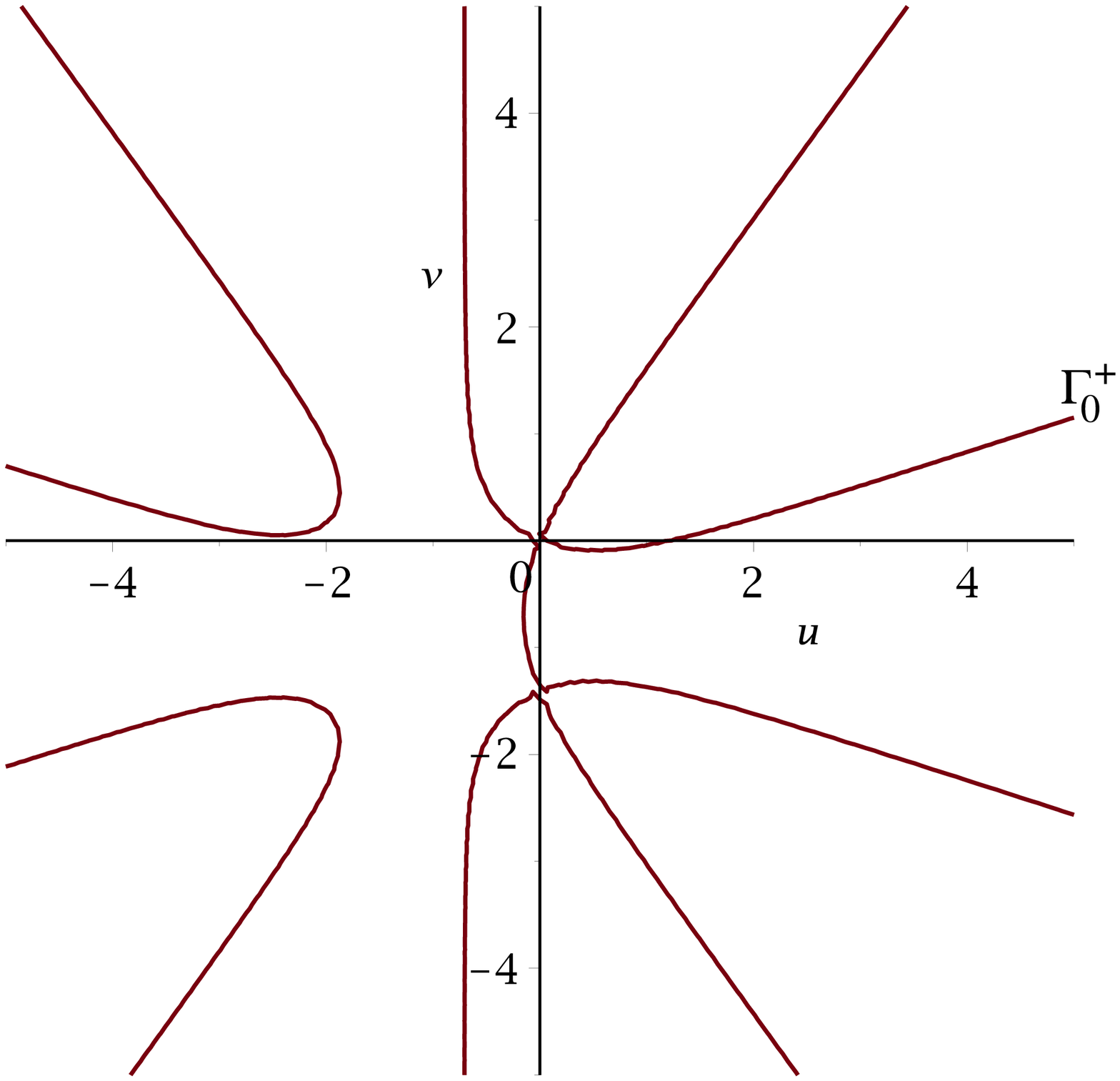}
\end{center}
\caption{\footnotesize{Curves $\textrm{Im}\, i(f(t) - f(t_k)) = 0$ for $k = 1$
(left) and $k = 0$ (right). In each figure, the coordinate used is $t - t_k = u + iv$, $u$, $v$ real, so that the saddle point at the origin in each plot.}}
\end{figure}

When $z < 0$, the steepest descent curve through $t_2$, $\Gamma_2^{-}$, is the one beginning at $\infty e^{9\pi i/10}$ passing through $t_2$ and dropping into the lower half plane to end at $\infty e^{13\pi i/10}$. The steepest descent curve through $t_3$, $\Gamma_3^{-}$, is the one beginning at $\infty e^{13\pi i/10}$ rising up to $t_3$ and then dropping again into the lower half-plane to end at $\infty e^{17\pi i/10}$. The steepest curve passing through $t_0$ (cf., the bottom illustration in Fig.~2), $\Gamma_0^{-}$, begins at $\infty e^{17\pi i/10}$, rises up to $t_0$ and then continues into the right of the upper half-plane to end at $\infty e^{\pi i/10}$. So, for $S(0,0,z)$ with $z < 0$, the original integration contour gets deformed into a sum of three steepest descent curves, $\Gamma_2^{-}+\Gamma_3^{-}+\Gamma_0^{-}$, in turn those passing through $t_2$, $t_3$ and $t_0$, yielding three relevant saddle points for our asymptotic analysis.

With the steepest descent paths and relevant saddles identified, we can construct the asymptotics of $Q$ for large $|z|$.

\begin{figure}[h]%
\begin{center}
\leavevmode
\epsfxsize=2in\epsfbox{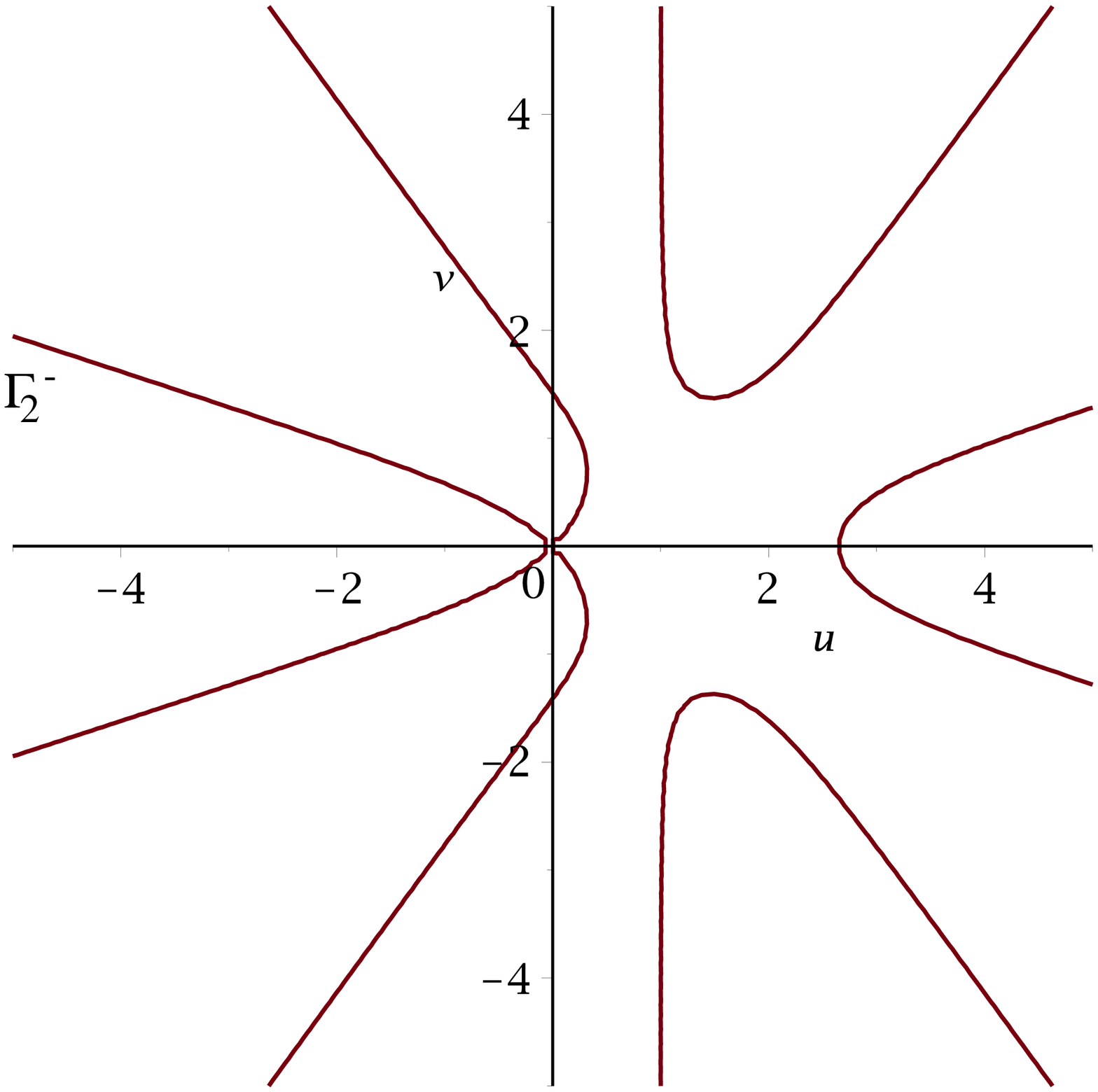}\quad\epsfxsize=2in\epsfbox{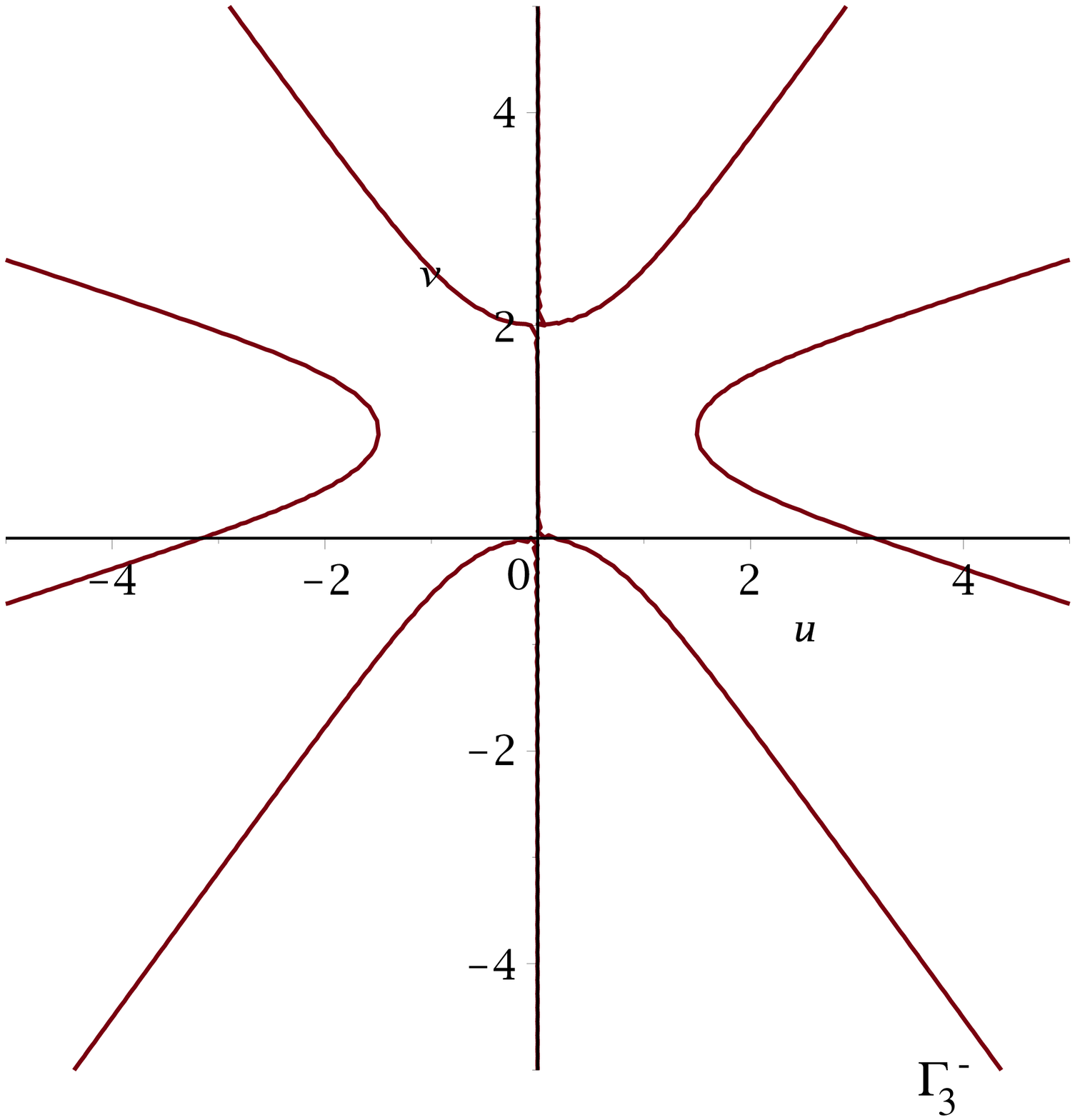}
\quad\epsfxsize=2in\epsfbox{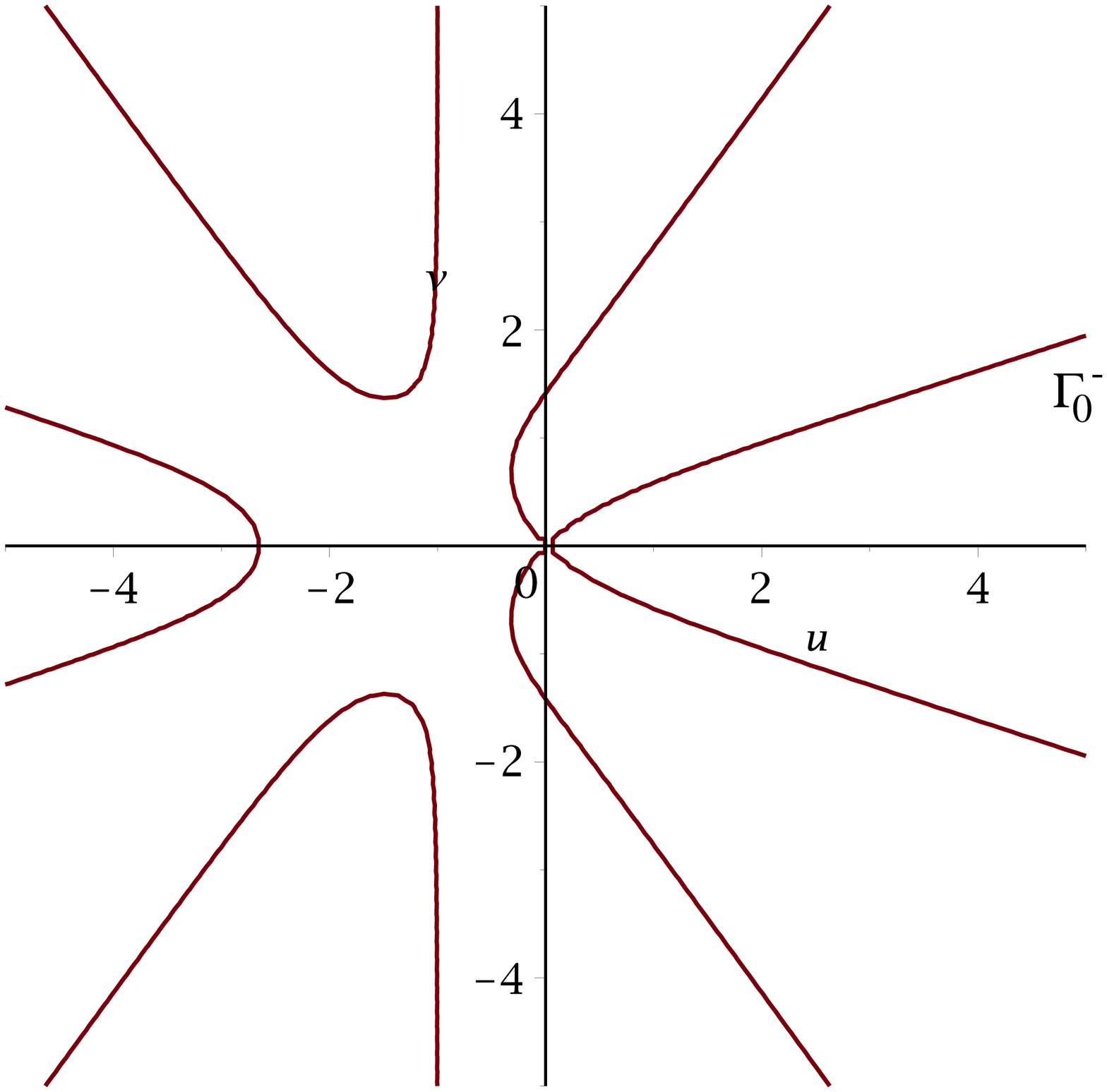}
\end{center}
\caption{\footnotesize{Curves $\textrm{Im}\, i(f(t) - f(t_k)) = 0$ for $k = 2$
(left) and $k = 3$ (right) and $k = 0$ (bottom). In each figure, the coordinate used is $t - t_k = u + iv$, $u$, $v$ real, so that the saddle point at the origin in each plot.}}
\end{figure}

For $z > 0$, there are only two relevant saddles. Writing $\lambda = z^{5/4}$, our steepest descent calculation proceeds along the usual lines (cf.~\cite[Ch.~II, {\S}4]{W}:
\begin{displaymath}
   \int_{\Gamma_k^{+}} e^{i\lambda (f_{+}(t) - f_{+}(t_k) + f_{+}(t_k))}dt
   = e^{i\lambda f_{+}(t_k)}\int_{\Gamma_k^{+}} e^{i\lambda(f_{+}(t) -
   f_{+}(t_k))}dt, k = 0, 1.
\end{displaymath}
With $i(f_{+}(t) - f_{+}(t_k)) = -\tau$ we find $t - t_k \sim \pm\sqrt{2i/f_{+}''(t_k)}\cdot \tau^{1/2}$ to leading order, and so 
\begin{displaymath}
   \frac{dt^{+}}{d\tau} - \frac{dt^{-}}{d\tau} \sim \sqrt{\frac{2i}{f_{+}''(t_k)}} \tau^{-1/2}
\end{displaymath}
to leading order. Accordingly,
\begin{displaymath}
   \int_{\Gamma_1^{+}+\Gamma_0^{+}} e^{i\lambda f_{+}(t)} dt
   \sim \left\{e^{i\lambda f_{+}(t_0)}\sqrt{\frac{2i}{4t_0^3}} 
      + e^{i\lambda f_{+}(t_1)}\sqrt{\frac{2i}{4t_1^3}}\right\}\sqrt{\frac{\pi}{\lambda}}.
\end{displaymath}
Since $2i/(4t_0^3) = \f{1}{2}e^{-\pi i/4}$ and $2i/(4t_1^3) = \f{1}{2}e^{\pi i/4}$,
\begin{displaymath}
   e^{i\lambda f_{+}(t_0)}\sqrt{\frac{2i}{4t_0^3}} 
      + e^{i\lambda f_{+}(t_1)}\sqrt{\frac{2i}{4t_1^3}}
   = \frac{e^{-4\lambda/5\sqrt{2}}}{\sqrt{2}}(e^{4\lambda i/(5\sqrt{2})-\pi i/8}
   + e^{-4\lambda i/(5\sqrt{2})+\pi i/8})
\end{displaymath}
we obtain
\begin{displaymath}
\int_{\Gamma_1^{+}+\Gamma_0^{+}} e^{i\lambda f_{+}(t)} dt
   \sim \sqrt{\frac{2\pi}{\lambda}} e^{-4\lambda/5\sqrt{2}}\cos(\frac{4\lambda}{5\sqrt{2}} - \frac{\pi}{8})
\end{displaymath}
as $\lambda \rightarrow \infty$, to leading order. Since
\begin{displaymath}
   Q(0,0,z) = z^{1/4} \int_{\Gamma_1^{+}+\Gamma_0^{+}} e^{i\lambda f_{+}(t)} dt
\end{displaymath}
if we want the asymptotic distribution of the zeroes of $Q$, then we will
want the cosine factor to vanish, so that we must have
\begin{displaymath}
   \frac{4\lambda}{5\sqrt{2}} - \frac{\pi}{8} = (2m+1)\frac{\pi}{2}
\end{displaymath}
for $m = 0, 1, 2, \ldots$. Isolating $\lambda$ and restoring the large
parameter $z$ results in the approximate zeroes
\begin{equation}\label{poszzeroes}
   z \sim \left[\frac{5\sqrt{2}}{4}\left(\frac{\pi}{8} + (2m+1)\frac{\pi}{2}\right)\right]^{4/5}, \quad m = 0, 1, 2, \ldots.
\end{equation}

For $z < 0$, three saddles are relevant and in an entirely similar manner as that for positive $z$, we have, with $\lambda = |z|^{5/4}$,
\begin{displaymath}
   \int_{\Gamma_k^{-}}e^{i\lambda f_{-}(t)}dt = 
   e^{i\lambda f_{-}(t_k)}\int_{\Gamma_k^{-}} e^{i\lambda (f_{-}(t) -
   f_{-}(t_k))} dt, \quad k = 2, 3, 0.
\end{displaymath}
Once again we are setting $i(f_{-}(t) - f_{-}(t_k)) = -\tau$ and to leading order we have $t - t_k \sim \pm\sqrt{2i/f_{-}''(t_k)} \tau^{1/2}$. We find
\begin{displaymath}
   \int_{\Gamma_2^{-}+\Gamma_3^{-}+\Gamma_0^{-}}e^{i\lambda f_{-}(t)}dt
   \sim
   \left\{ e^{i\lambda f_{-}(t_0)}\sqrt{\frac{i}{2}} + e^{i\lambda f_{-}(t_2)}\sqrt{\frac{i}{-2}} + e^{i\lambda f_{-}(t_3)}\sqrt{\frac{i}{2i}}\right\}\sqrt{\frac{\pi}{\lambda}}
\end{displaymath}
since the expressions for $t_k$ are so much simpler in the $z < 0$ case.
With the principal branch of square root being used, we find
\begin{displaymath}
   \int_{\Gamma_2^{-}+\Gamma_3^{-}+\Gamma_0^{-}}e^{i\lambda f_{-}(t)}dt
   \sim \sqrt{\frac{\pi}{2\lambda}}\left\{ e^{-4\lambda i/5 + \pi i/4}
   + e^{4\lambda i/5 - \pi i/4} + e^{-4\lambda/5}\right\}
\end{displaymath}
and since the last term is exponentially negligible compared to the other two, we arrive at the first order approximation
\begin{displaymath}
   \int_{\Gamma_2^{-}+\Gamma_3^{-}+\Gamma_0^{-}}e^{i\lambda f_{-}(t)}dt
   \sim \sqrt{\frac{2\pi}{\lambda}} \cos\left(\frac{4\lambda}{5} -
   \frac{\pi}{4}\right)
\end{displaymath}
as $\lambda \rightarrow +\infty$. If this is to vanish, then we must have
\begin{displaymath}
   \frac{4\lambda}{5} - \frac{\pi}{4} = (2m+1)\frac{\pi}{2}
\end{displaymath}
for integers $m$, and restoring $z$,
\begin{equation}\label{negzzeroes}
   z \sim -\left[\frac{5}{4}\left(\frac{\pi}{4} + (2m+1)\frac{\pi}{2}\right)
   \right]^{4/5}, \quad m = 0, 1, 2, \ldots.
\end{equation}

The approximate zeroes in (\ref{poszzeroes}) and (\ref{negzzeroes}) are the zeroes identified as the line of zeroes in \cite{CCF}.

\section{$Q(0,y,z)$}
In allowing $y$ to be nonzero, we begin as before by rescaling the integration variable with $t \mapsto |z|^{1/4}t$ to produce
\begin{displaymath}
   Q(0,y,\pm|z|) = |z|^{1/4}\int_{\infty e^{9\pi i/10}}^{\infty e^{\pi i/10}} 
      e^{i|z|^{5/4}f_{\pm}(t; \gamma)}dt
\end{displaymath}
where now
\begin{equation}\label{fgamma}
   f_{\pm}(t; \gamma) = f_{\pm}(t) = \f{1}{5}t^5 + \f{1}{2}\gamma t^2 \pm t
\end{equation}
and
\begin{displaymath}
   \gamma = \frac{y}{|z|^{3/4}}.
\end{displaymath}
If it happens that we restrict $y$ to be bounded, then the quantity $\gamma$ is clearly evanescent as $|z| \rightarrow \infty$, and so the study of $Q(0,y,z)$ reduces to that of $Q(0,0,z)$ for which there are, indeed only zeroes along the $z$-axis.

A more subtle approach is needed if $\gamma > 0$; recall (\ref{symmetry}).
If $(y,z)$ lies above the caustic in the $yz$-plane, $(z/3)^3 = (y/4)^4$,
corresponding to $\gamma = 4/3^{3/4}$, then we know the phase function
(\ref{fgamma}) $f_{+}(t)$ has saddle points consisting of two complex conjugate pairs, say 
\begin{displaymath}
   t_1 = p_1 + iq_1, \quad t_2 = p_1 - iq_1 \quad \textrm{and} \quad
   t_3 = p_2 + iq_2, \quad t_4 = p_2 - iq_2,
\end{displaymath}
where $p_1, p_2, q_1, q_2$ are real numbers with $q_1$ and $q_2$ nonnegative.
From the theory of equations, we know that the sum of these saddle points must equal the coefficient of $t^3$ in $f_{+}'(t)$, or
$t_1 + t_2 + t_3 + t_4 = 0$. Thus, $2p_1 + 2p_2 = 0$ whence $p_1 = -p_2$. Let us take $p = p_1 \geq 0$ and so $p_2 = -p$ and our roots of $f_{+}' = 0$ have the form
\begin{equation}\label{roots}
   t_1 = p + iq_1, \quad t_2 = p - iq_1, \quad t_3 = -p + iq_2,
   \quad t_4 = -p - iq_2
\end{equation}
with all of $p$, $q_1$, $q_2$ nonnegative.

Additionally, the coefficient of $t^2$ in $f_{+}' = 0$ is also 0, so the elementary symmetric function
\begin{displaymath}
   t_1t_2 + t_1t_3 + t_1t_4 + t_2t_3 + t_2t_4 + t_3t_4 
\end{displaymath} 
must vanish. Using the values (\ref{roots}) in this symmetric function leads 
to the identity $q_1^2 + q_2^2 = 2p^2$.

In considering the steepest descent curve structure applicable in this case, we see that the circumstance for $z > 0$ in our analysis of $Q(0,0,z)$ applies directly in this case: the relevant saddle points contributing to the asymptotic behaviour of $Q(0,y,z)$ are the two saddles above the real axis in the $t$-plane, so we need only consider $t_1$ and $t_3$. Furthermore, the steepest descent paths then have the same form as $\Gamma_1^{+}$ and $\Gamma_0^{+}$ of {\S}2 and in Fig.~1. Therefore, the dominant contributions to the asymptotic behaviour of $Q(0,y,z)$  are of the form
\begin{equation}\label{contribution}
   e^{i\lambda f_{+}(t_k)+\pi i/4 - i(\arg f_{+}''(t_k))/2}
   \sqrt{\frac{2\pi}{\lambda |4t_k^3 + \gamma|}}
\end{equation}
for $k = 1, 3$ with $\lambda = z^{5/4}$; here, $f_{+}''(t_k) = 4t_k^3 + \gamma$. Since $t_k^4 + \gamma t_k + 1 = 0$, we find that $t_k^5 = -t_k-\gamma t_k^2$ and so $f_{+}(t_k) = \f{3}{10}\gamma t_k^2 + \f{4}{5}t_k$ and evaluating at $t_1$ and $t_3$ gives
\begin{eqnarray*}
   f_{+}(t_1) & = & \{\f{3}{10}\gamma(p^2-q_1^2)+\f{4}{5}p\} + i\{\f{3}{10}\gamma\cdot 2pq_1 + \f{4}{5}q_1\} \\
   f_{+}(t_3) & = & \{\f{3}{10}\gamma(p^2-q_2^2)-\f{4}{5}p\} - i\{\f{3}{10}\gamma\cdot 2pq_2 - \f{4}{5}q_1\}.
\end{eqnarray*}
In light of the relation $q_1^2 + q_2^2 = 2p^2$, we have $p^2 - q_2^2 =
-(p^2 - q_1^2)$ so that these evaluations lead to
\begin{eqnarray*}
   if_{+}(t_1) & = & -\{\f{3}{10}\gamma\cdot 2pq_1 + \f{4}{5}q_1\} + i\{\f{3}{10}\gamma(p^2 - q_1^2) + \f{4}{5}p\}\\
   if_{+}(t_3) & = & +\{\f{3}{10}\gamma\cdot 2pq_2 - \f{4}{5}q_2\} - i\{\f{3}{10}\gamma(p^2 - q_1^2) + \f{4}{5}p\}.
\end{eqnarray*}
If these are used in (\ref{contribution}), then we see that the contribution of the saddle point $t_1$ is of exponentially small order compared to the contribution from $t_3$ and so there is no way to combine the two contributions in a form that would permit us to extract a zero of $Q(0,y,z)$,
unless we had $p = 0$.

If it were to happen that $p = 0$, then there would be four saddle points for our integral strung along the imaginary axis in the $t$-plane, an impossible occurrence under the current hypotheses, for that then implies that $f_{+}'(t)$ has a nonzero quadratic term which can only happen in the case where $x \neq 0$ in $Q(x,y,z)$.

So, above the caustic, it appears the only zeroes of $Q(0,y,z)$ for $z > 0$ occur along the $z$-axis.

For $(y,z)$ below the caustic, the saddle point arrangement changes. To fix our discussion, we assume $z < 0$ (which is certainly below the caustic) and the relevant phase function now is $f_{-}(t;\gamma) = f_{-}(t)$. As we saw in {\S}2, in this setting, $f_{-}'(t) = 0$ has one pair of real roots, and one complex conjugate pair. Let the real roots be $t_1 < t_2$, and let $t_3$ and $t_4$ be the conjugate pair, say $t_3 = p + iq$ and $t_4 = p - iq$ with $q > 0$. The arrangement of saddle points and steepest descent curves we saw in {\S}2 for $z < 0$ carries over to this case with $y > 0$: the integration contour for $Q(0,y,z)$ can be deformed into a sum of steepest paths, one through $t_1$, one through $t_4$ and one through $t_2$, as was the case for $\Gamma_2^{-}+\Gamma_3^{-}+\Gamma_0^{-}$; recall Fig.~2.

From the theory of equations, we know $\sum t_k = 0$ whence $t_1 + t_2 = -2p$
and $\sum_{j<k}t_jt_k = 0$ implies $t_1t_2 = 3p^2 - q^2$.

Since $t_1$ and $t_2$ are real saddles, $f_{-}$ must have a local max at $t_1$ and a local min at $t_2$, and so $f_{-}''(t_1) = 4t_1^3 + \gamma < 0$ and $f_{-}''(t_2) = 4t_2^3 + \gamma > 0$. As well, $f_{-}'(t_k) = 0$ implies that $t_k^4 = 1 - \gamma t_k$ giving the evaluations $f_{-}(t_k) = \f{3}{10}\gamma t_k^2 - \f{4}{5}t_k$. The saddle at $t_4$ will result in an exponentially negligible contribution, and so we find, to leading order, that
the saddles $t_1$ and $t_2$ contribute
\begin{displaymath}
   e^{i\lambda f_{-}(t_1)-\pi i/4}\sqrt{\frac{2\pi}{\lambda|f_{-}''(t_1)|}}
+  e^{i\lambda f_{-}(t_2)+\pi i/4}\sqrt{\frac{2\pi}{\lambda f_{-}''(t_2)}}
\end{displaymath}
where $\lambda = |z|^{5/4}$, as before.

If these contributions are to combine into a cosine term as in the previous cases we've examined, then we will need to have
\begin{displaymath}
   f_{-}(t_1) = -f_{-}(t_2) \quad \textrm{and} \quad
   4t_1^3 + \gamma = -(4t_2^3 + \gamma).
\end{displaymath}
The first of this pair then implies
\begin{displaymath}
   \f{3}{10}\gamma t_1^2 - \f{4}{5}t_1 = -(\f{3}{10}\gamma t_2^2 - \f{4}{5}t_2)
\end{displaymath}
or
\begin{displaymath}
   0 = \f{3}{10}\gamma(t_1^2 + t_2^2) + \f{8}{5}p.
\end{displaymath}
Since $\gamma > 0$ and $t_1$ and $t_2$ are real, then we must have $p < 0$.
$t_1 + t_2 = -2p$ then implies that $t_2 > |t_1|$. To have 
$4t_1^3 + \gamma = -(4t_2^3 + \gamma)$ is equivalent to $-2(t_1^3 + t_2^3) =
\gamma$. But $t_2 > |t_1|$ gives $t_2^3 > |t_1|^3$ from which we obtain
$t_2^3 + t_1^3 > 0$ which in turn implies that $\gamma < 0$, a contradiction.

Therefore the contributions of the saddle points $t_1$ and $t_2$ to the asymptotics of $Q(0,y,z)$ for $z < 0$ cannot add in a way as to produce a cosine factor.

Therefore, the only zeroes of $Q(0,y,z)$ for $z$ of either sign lie along the $z$-axis.

\section{Closing remarks}
Pearcey \& Hill \cite[p.~48]{PH} assert that if
\begin{displaymath}
   I_5(X,Y) = \int_{-\infty}^\infty e^{i(t^5 + Xt^2 + Yt)} dt,
\end{displaymath}
then $I_5(0,Y)$ has zeroes on the positive $Y$-axis near
\begin{displaymath}
   Y= 5\cdot 2^{-6/5}(n+\f{5}{8})^{4/5}\cdot \pi^{4/5},
\end{displaymath}
for $n = 0, 1, 2, \ldots$, a result that we have recovered in our (\ref{poszzeroes}), once the change in parameters in (\ref{SQ}) and (\ref{QS}) is taken into account. However, our result for the zeroes of $Q(0,0,z)$, given in (\ref{negzzeroes}) for negative $z$, although close to what is reported in \cite[p.~52]{PH}, suggests an arithmetic error in \cite{PH} -- there is an incorrect scaling factor there.

The analysis of the distribution of zeroes of the swallowtail integral for $x \neq 0$ is a more complicated affair, which we elect to leave for another time.

Finally, we close by noting the relevence of \cite{Nye} to this discussion. Nye discusses families of zeroes of the swallowtail integral, using the location of saddlepoints of $S(x,0,z)$ to anchor his analysis. However, his analysis does not appear to produce explicit formul\ae\  for the families of zeroes of $S(x,y,z)$, though he is able to provide a means of indexing families of these zeroes.


\begin{thebibliography}{NIST}
\bibitem[BK]{BK}{\sc M V Berry and S Klein}, {\em Colored diffraction
catastrophes,\/} Proc Natl Acad Sci USA, Vol~93 (1996), pp~2614-2619

\bibitem[CCF]{CCF}{\sc J N L Connor, P R Curtis and D Farrelly}, {\em The
uniform asymptotic swallowtail approximation: practical methods for oscillating
integrals with four coalescing saddle points,\/} J Phys A: Math Gen, Vol~17
(1984), pp~283-310

\bibitem[CH]{CH}{\sc J N L Connor and C A Hobbs}, {\em Numerical evaluation of cuspoid and Bessoid oscillating integrals for applications in chemical physics,\/} {\cyr Khimicheskaya Fizika}, Vol~23 (2004), No~2, pp~13-19

\bibitem[Kam]{Kam}{\sc D Kaminski}, {\em Asymptotic expansion of the
Pearcey integral near the caustic,\/} SIAM J Math Anal, Vol~20 (1989), No~4, 
pp~987-1005

\bibitem[NIST]{NIST} {\sc F W J Olver, D W Lozier, R F Boisvert and C W Clark (eds.)}, {\em The NIST Handbook of Mathematical Functions,\/}
National Institute of Standards and Technology, US Department of Commerce, and Cambridge University Press, Cambridge, 2010

\bibitem[Nye]{Nye}{\sc J F Nye}, {\em Dislocation lines in the swallowtail diffraction catastrophe,\/} Proc Roy Soc A (2007), Vol 463, pp~343-355

\bibitem[PH]{PH}{\sc T Pearcey and G W Hill}, {\em The Effect of Aberrations upon the Optical Focus,\/} Commonwealth Scientific and Industrial Research Organization, 1963

\bibitem[W]{W}{\sc R Wong}, {\em Asymptotic Approximations of Integrals,
\/} Academic Press, Toronto, 1989
\end{thebibliography}
\end{document}